\documentclass[10pt]{article}
\usepackage{amsmath,amssymb,amsthm,amscd}
\numberwithin{equation}{section}

\def\p{\partial}
\def\b{\bar}

\def\o{\omega}
\def\l{\lambda}

\def\cL{{\cal L}}

\def\cC{{\cal C}}

\def\cP{{\cal P}}

\def\cA{{\mathcal A}}
\def\cC{{\mathcal C}}

\def\cL{{\mathcal L}}

\def\cN{{\mathcal N}}

\def\cP{{\mathcal P}}

\def\cR{{\mathcal R}}

\def\RR{{\mathbb R}}

\newtheorem{prop}{Proposition}
\newtheorem{theo}[prop]{Theorem}
\newtheorem{lem}[prop]{Lemma}
\newtheorem{cor}[prop]{Corollary}
\newtheorem{rem}[prop]{Remark}

\newtheorem{defi}[prop]{Definition}

\newtheorem{q}[prop]{Question}
\def\begeq{\begin{equation}}
\def\endeq{\end{equation}}

\def\and{\quad{\rm and}\quad}

\let\lra=\longrightarrow

\def\mapright\#1{\,\smash{\mathop{\lra}\limits^{\#1}}\,}

\begin {document}
\bibliographystyle{plain}
\title{On the lower bound of energy functional $E_1$ (I)\\---a stability theorem on the K\"ahler Ricci flow}
\author{Xiuxiong  Chen}

\date{
}
\maketitle
\begin{center} {\bf In memory of the great mathematician S. S. Chern}
\end{center}
\tableofcontents
\newpage
\section{Introduction and main results}
This is the first of a serious papers aiming to study the lower
bound of the energy functional $E_1$\footnote{ This is first
introduced in \cite{chentian001}. For the convenience of readers,
we will give precise definition of $E_1$ in Section 2 (cf.
Definition \ref{defi:E1} below).} and its relation to the
convergence of the K\"ahler Ricci flow as well as the relation to
the existence of K\"ahler Einstein metrics. The
K\"ahler Ricci flow is known for exists for all times and
converges to K\"ahler Einstein metrics when the first Chern class
was negative or zero \cite{Cao85} \cite{cao92}.  When the first
Chern class is positive, there is very few known results on global
convergence.   In the case of positive bisectional curvature, the
convergence of the flow  is first proved in \cite{chentian001}
\cite{chentian002} where the energy functional $E_1$ plays a key role. After the important work of Perelman, there
are active researches in this direction and many good works
appear, for instance, \cite{Cao03} \cite{Natasa2} (more recently,
\cite{phongjacob042})and references therein. In the present paper, we
 prove a stability theorem of the K\"ahler Ricci flow near
the infimum of $E_1$ under the assumption  that the initial
metric has  Ricci $> -1$  and $|Riem|$ bounded.
 The underlying moral is: if a K\"ahler metric is
sufficiently closed
to a K\"ahler Einstein metric,  then the K\"ahler Ricci flow shall converges to it.
The present work should be viewed as first step in a more ambitious program of deriving the existence
of K\"ahler Einstein metric with arbitrary energy level, provided that this energy functional
has a uniform lower bound in this K\"ahler class.  \\

  This stability theorem can be alternatively viewed as generalization of the
well known  pointwise  curvature-pinching theorem proved by method of Ricci
flow in the  1980s. The beautiful work  of R. Hamilton \cite{Hamilton82}  proved by Ricci flow, that any closed 3-manifold of positive Ricci curvature is
diffeomorphic to a spherical space form.  This result was
generalized by several mathematicians (c.f., \cite{Hu85},
\cite{Margein86} \cite{Nish86} ) to higher dimensional manifold:
if the Riemannian curvature satisfies certain pointwise pinching
condition, then the flow will preserve this pinching condition and
converges to a spherical space form up re-scaling of the evolved
metrics \footnote{In dimension 4, there is a sharp pinching
condition due to C. Margerin \cite{Margein98}.}. One wonders if  a similar result in
the K\"ahler setting is feasible. Unfortunately, the corresponding pinching
condition doesn't hold in K\"ahler setting in general. The present
work can be viewed as a generalization in this direction with one
crucial difference:  We didn't prove that the pinching condition
imposed initially is preserved over time. Instead,  we prove that,
even though the pinching condition might lost immediately after the flow
starts, it will be recovered after a fixed period of time (the
length of this period is determined by the geometrical condition of the initial metric.). This periodic re-visit of the pinching condition
 is enough to give us control of  $L^\infty$ norm of the curvature of the evolving metrics.  In many
cases, this will lead to the convergence of the K\"ahler Ricci flow.\\

More specifically, let $(M,\omega_0)$ be a polarized $n-$ dimensional compact K\"ahler
manifold with positive first Chern class where $[\omega_0]$ is the
canonical K\"aher class.  We always normalize the K\"ahler class
so that the canonical class is the Ricci class. One of the main theorems of
this note is to prove the following stability theorem:
\begin{theo} \label{th:main1} For any $\delta, \Lambda > 0, $ there exists a small positive constant $\epsilon (\delta, \Lambda) > 0$
 such that if the subspace $\cA(\delta, \Lambda, \epsilon)$ of K\"ahler metrics in $[\omega]$
\begin{equation}
 \{\omega_g \in [\omega_0]\mid
Ric(g)
> -1 +\delta, \; \mid Riem(g) \mid < \Lambda, \; E_1(g)
\leq \displaystyle \inf_{Ric(\omega)\geq -1}\; E_1(\omega)
+\epsilon\} \label{eq:stabilitycondition1}
\end{equation}
is non-empty, then there exists a K\"ahler Einstein metrics in
this canonical K\"ahler class if the following integral condition
hold
\begin{equation}
[C_1(M)]^{[2]}\cdot [\omega]^{[n-2]} - {{2(n+1)}\over n} [ C_2(M)]
\cdot [\omega]^{[n-1]} = 0. \label{intro:topologycondition0}
\end{equation}
 Moreover, for any metric $g_1\in \cA(\delta, \Lambda,
\epsilon)$, the K\"ahle Ricci flow will deform it exponentially
fast to a K\"ahler-Einstein metric in the limit.
\end{theo}
\begin{rem} The reason to assume a $L^\infty$ bound in the Riemannian curvature of the initial
metric is to ensure the inequality $Ric>-1$ is preserved for some
fixed amount of time.  One should be able to replace the
$L^\infty$ bound on the bisectional curvature by some $L^p$
estimates for some $p > {{n+2}\over 2}.\;$ Condition \ref{intro:topologycondition0} is introduced to control the
norm of sectional curvature. This condition is removed in Theorem \ref{th:2dim} below.
\end{rem}
\begin{rem} In Theorem \ref{th:main1},
we can estimate $\epsilon(\Lambda, \delta)$ explicitly:
\begin{equation}
\epsilon(\delta, \Lambda)  \leq \left({1\over \Lambda}
\right)^{2n} \delta \cdot\epsilon_0(n)^2. \label{eq:bestconstant1}
\end{equation}
 Here
$\epsilon_0(n)$ is some universal constant which depends only on
the dimension (cf. Lemma \ref{th:bestconstant2}). This estimate
is not optimal.
\end{rem}

This naturally leads to the following  question.
\begin{q} In any canonical K\"ahler class $\omega_0$ where the first
Chern class is positive, we can define the following
invariant
\[
\beta(M, \omega) = \displaystyle \inf_{\varphi\in \cP(M, [\omega])} \;\displaystyle \max_{x\in M} |Ric(\omega_\varphi) -\omega_\varphi|_{\varphi}.
\]
Then this defined a holomorphic invariant which depends only on the underlying complex structure.   Apparently, if there is K\"ahler
Einstein metric in $[\omega_0]$, then this invariant vanishes.
What about the manifold without K\"ahler Einstein metric but this invariant vanishes? How does this invariant change
with respect to the deformation of complex structure?  An obvious guess is that  this is related to the stability of
the tangent bundle.   When the first Chern class is negative, Aubin, Yau's solution of Calabi conjecture implies
that this invariant (defined accordingly) always vanishes. 
\end{q}

\begin{theo} \label{th:general} Conditions as stated in Theorem \ref{th:main1} except
that the assumption \ref{intro:topologycondition0} is replaced by
the following: there exists a constant $C(p)$ such that
\begin{equation} \displaystyle \int_M\; |Riem(g(t))|_{g(t)}^p
\;d\,vol_g \leq C, \qquad {\rm for\; fixed \;} p > n+1.  \label{eq:lpcurvatureassumption}
\end{equation}
Here $g(t)$ is the evolved K\"ahler metrics along the K\"ahler Ricci
flow initiated from $g(0).\;$  Then the K\"ahler Ricci flow is
non-singular (i.e.,, the bisectional curvature of $g(t)$ are
uniformly bounded).  Moreover,  the $Ric(g(t))-g(t)$ converges to
$0$ uniformly as $t \rightarrow \infty$.  If in additional, we
assume that the underlying complex structure is stable in the sense
that no complex structure in the closure of its orbit of
diffeomorphism group contains larger holomorphic automorphism group,
then the flow converges to a  unique K\"ahler Einstein metric in the
original complex structure exponentially fast.
\end{theo}

A natural question is whether condition
\ref{eq:lpcurvatureassumption} will actually occur?
\begin{theo} \label{th:2dim} Conditions as in Theorem \ref{th:main1} except the
assumption \ref{intro:topologycondition0} is removed.  In
additional, we assume $M$ is complex surface and $X \neq 0$ is a
holomorphic vector field in $M$ such that $L_{Im (X)}\; g_0 =
0.\;$ Then, the K\"ahler Ricci flow
converges to a K\"ahler Einstein metrics in the limit
exponentially fast.
\end{theo}

\begin{rem} Tian's solution \cite{Tian87}  of Calabi conjecture in complex surface with positive first Chern class certainly
include this as a special case. This might be interesting since it demonstrates how to obtain inequality \ref{eq:lpcurvatureassumption}
in absence of topological condition \ref{intro:topologycondition0}.
\end{rem}
 {\bf
Acknowledgment}  Main results of this paper has been reported
in Nan Jing University and Tokyo Institute of Technology,  during
May-June, 2004. The author wishes to thank both University for their
hospitality.  The author also wishes to thank both Professor S. K. Donaldson and G. Tian for their
 constant, warm encouragements on my research in last few years.
\section{ Basic  K\"ahler geometry}
\subsection{Setup of notations}

Let $M$ be an $n$-dimensional compact K\"ahler manifold. A K\"ahler metric can be given by its
K\"ahler form $\omega$ on $M$. In local coordinates $z_1, \cdots, z_n$, this $\omega$ is of
the form
\[
\omega = \sqrt{-1} \displaystyle \sum_{i,j=1}^n\;g_{i
\overline{j}} d\,z^i\wedge d\,z^{\overline{j}}  > 0,
\]
where $\{g_{i\overline {j}}\}$ is a positive definite Hermitian matrix function.
The K\"ahler condition requires that $\omega$ is a closed positive
(1,1)-form. In other words, the following holds
\[
 {{\partial g_{i \overline{k}}} \over
{\partial z^{j}}} =  {{\partial g_{j \overline{k}}} \over
{\partial z^{i}}}\qquad {\rm and}\qquad {{\partial g_{k
\overline{i}}} \over {\partial z^{\overline{j}}}} = {{\partial
g_{k \overline{j}}} \over {\partial
z^{\overline{i}}}}\qquad\forall\;i,j,k=1,2,\cdots, n.
\]
The K\"ahler metric corresponding to $\omega$ is given by
\[
 \sqrt{-1} \;\displaystyle \sum_1^n \; {g}_{\alpha \overline{\beta}} \;
d\,z^{\alpha}\;\otimes d\, z^{ \overline{\beta}}.
\]
For simplicity, in the following, we will often denote by $\omega$ the corresponding K\"ahler metric.
The K\"ahler class of $\omega$ is its cohomology class $[\omega]$ in $H^2(M,\RR).\;$
By the Hodge theorem, any other K\"ahler
metric in the same K\"ahler class is of the form
\[
\omega_{\varphi} = \omega + \sqrt{-1} \displaystyle \sum_{i,j=1}^n\;
{{\partial^2 \varphi}\over {\partial z^i \partial z^{\overline{j}}}}
> 0
\]
for some real valued function $\varphi$ on $M.\;$ The functional
space in which we are interested (often referred as the space of
K\"ahler potentials) is
\[
{\cal P}(M,\omega) = \{ \varphi \;\mid\; \omega_{\varphi} = \omega
+ \sqrt{-1}
 {\partial} \overline{\partial} \varphi > 0\;\;{\rm on}\; M\}.
\]
Given a K\"ahler metric $\omega$, its volume form  is
\[
  \omega^n = {1\over {n!}}\;\left(\sqrt{-1} \right)^n \det\left(g_{i \overline{j}}\right)
 d\,z^1 \wedge d\,z^{\overline{1}}\wedge \cdots \wedge d\,z^n \wedge d\,z^{\overline{n}}.
\]
Its Christoffel symbols are given by
\[
  \Gamma^k_{i\,j} = \displaystyle \sum_{l=1}^n\;g^{k\overline{l}} {{\partial g_{i \overline{l}}} \over
{\partial z^{j}}} ~~~{\rm and}~~~ \Gamma^{\overline{k}}_{\overline{i}\,\overline{j}} =
\displaystyle \sum_{l=1}^n\;g^{\overline{k}l} {{\partial g_{l \overline{i}}} \over
{\partial z^{\overline{j}}}}, \qquad\forall\;i,j,k=1,2,\cdots n.
\]
The curvature tensor is
\[
 R_{i \overline{j} k \overline{l}} = - {{\partial^2 g_{i \overline{j}}} \over
{\partial z^{k} \partial z^{\overline{l}}}} + \displaystyle \sum_{p,q=1}^n g^{p\overline{q}}
{{\partial g_{i \overline{q}}} \over
{\partial z^{k}}}  {{\partial g_{p \overline{j}}} \over
{\partial z^{\overline{l}}}}, \qquad\forall\;i,j,k,l=1,2,\cdots n.
\]
We say that $\omega$ is of nonnegative bisectional curvature if
\[
 R_{i \overline{j} k \overline{l}} v^i v^{\overline{j}} w^k w^{\overline{l}}\geq 0
\]
for all non-zero vectors $v$ and $w$ in the holomorphic tangent
bundle of $M$. The bisectional curvature and the curvature tensor
can be mutually determined. The Ricci curvature of $\omega$ is
locally given by
\[
  R_{i \overline{j}} = - {{\partial}^2 \log \det (g_{k \overline{l}}) \over
{\partial z_i \partial \bar z_j }} .
\]
So its Ricci curvature form is
\[
  {\rm Ric}(\omega) = \sqrt{-1} \displaystyle \sum_{i,j=1}^n \;R_{i \overline{j}}(\omega)
d\,z^i\wedge d\,z^{\overline{j}} = -\sqrt{-1} \partial \overline{\partial} \log \;\det (g_{k \overline{l}}).
\]
It is a real, closed (1,1)-form. Recall that $[\omega]$ is called
a canonical K\"ahler class if this Ricci form is cohomologous to
$\lambda \;\omega,\; $ for some constant $\lambda.$ In our
setting, we require $\lambda = 1.\;$

 \subsection{The K\"ahler Ricci flow}

    Now we assume that the first Chern class $c_1(M)$ is positive.
The normalized Ricci flow (c.f. \cite{Hamilton82} and
\cite{Hamilton86}) on a K\"ahler manifold $M$ is of the form
\begin{equation}
  {{\partial g_{i \overline{j}}} \over {\partial t }} = g_{i \overline{j}}
  - R_{i \overline{j}}, \qquad\forall\; i,\; j= 1,2,\cdots ,n.
\label{eq:kahlerricciflow}
\end{equation}
If we choose the initial K\"ahler metric $\omega$ with $c_1(M)$ as
its K\"ahler class. The flow (2.1) preserves the K\"ahler class
$[\omega]$. It follows that on the level of K\"ahler potentials,
the Ricci flow becomes
\begin{equation}
   {{\partial \varphi} \over {\partial t }} =  \log {{\omega_{\varphi}}^n \over {\omega}^n } + \varphi - h_{\omega} ,
\label{eq:flowpotential}
\end{equation}
where $h_{\omega}$ is defined by
\[
  {\rm Ric}(\omega)- \omega = \sqrt{-1} \partial \overline{\partial} h_{\omega}, \; {\rm and}\;\displaystyle \int_M\;
  (e^{h_{\omega}} - 1)  {\omega}^n = 0.
\]
Then the evolution equation for bisectional curvature is

\begin{eqnarray}{{\partial }\over {\partial t}} R_{i \overline{j} k
\overline{l}} & = & \bigtriangleup R_{i \overline{j} k
\overline{l}} + R_{i \overline{j} p \overline{q}} R_{q
\overline{p} k \overline{l}} - R_{i \overline{p} k \overline{q}}
R_{p \overline{j} q \overline{l}} + R_{i
\overline{l} p \overline{q}} R_{q \overline{p} k \overline{j}} + R_{i \overline{j} k \overline{l}} \nonumber\\
& &  \;\;\; -{1\over 2} \left( R_{i \overline{p}}R_{p \overline{j}
k \overline{l}}  + R_{p \overline{j}}R_{i \overline{p} k
\overline{l}} + R_{k \overline{p}}R_{i \overline{j} p
\overline{l}} + R_{p \overline{l}}R_{i \overline{j} k
\overline{p}} \right). \label{eq:evolutio of curvature1}
\end{eqnarray}

The evolution equation for Ricci curvature and scalar curvature
are
 \begin{eqnarray} {{\p R_{i \b j}}\over {\p t}} & = & \triangle
  R_{i\b j} + R_{i\b j p \b q} R_{q \b p} -R_{i\b p} R_{p \b j},\label{eq:evolutio of curvature2}\\
  {{\p R}\over {\p t}} & = & \triangle R + R_{i\b j} R_{j\b i}- R.
  \label{eq:evolutio of curvature3}
  \end{eqnarray}

As usual, the flow equation (\ref{eq:kahlerricciflow}) or
(\ref{eq:flowpotential}) is referred as the K\"ahler Ricci flow on
$M$. It is proved by Cao \cite{Cao85}, who followed Yau's
celebrated work \cite{Yau78}, that the K\"ahler Ricci flow exists
globally for any smooth initial K\"ahler metric. It was proved by
S. Bando \cite{Bando84} in dimension 3 and N. Mok \cite{Mok88} in
all dimension that positivity of the bisectional curvature is
preserved under the flow. In \cite{chentian001}
\cite{chentian002}, we prove that the K\"ahler Ricci flow, in
K\"ahler Einstein manifold, initiated from a metric with positive
bisectional curvature converges to a K\"ahler Einstein metric with
constant bisectional curvature.

\subsection{The energy functional $E_1$} In
this subsection, we introduce the new functionals \cite{chentian001} $E_k =
E_k^0 - J_k (k =0,1,2\cdots, n) $ where $E_k^0$ and $J_k$ are
defined below.

\begin{defi} For any $k=0,1,\cdots, n$, we define a functional $E_k^0$
on ${\cal P}(M,\omega)$ by
\[
  E_{k,\omega}^0 (\varphi) = {1\over V}\; \displaystyle \int_M\;  \left( \log {{\omega_{\varphi}}^n \over \omega^n}
   - h_{\omega}\right) \left(\displaystyle \sum_{i=0}^k\; {{\rm Ric}(\omega_{\varphi})}^{i}\wedge\omega^{k-i} \right)
    \wedge {\omega_{\varphi}}^{n-k} + c_k,
\]
where
\[
c_k ={1\over V}\; \displaystyle \int_M\;
   h_{\omega} \left(\displaystyle \sum_{i=0}^k\; {{\rm Ric}(\omega)}^{i}\wedge\omega^{k-i} \right)
    \wedge {\omega}^{n-k}.
\]
\end{defi}

\begin{defi}\label{defi:E1}
For each $k=0,1,2,\cdots, n-1$, we will define $J_{k, \omega}$ as
follows: Let $\varphi(t) $ ($t\in [0,1]$) be a path from $0$ to
$\varphi$ in ${\cal P}(M,\omega)$, we define
\[ J_{k,\o}(\varphi) = -{n-k\over V}
\int_0^1 \int_M {{\partial \varphi}\over{\partial t}}
\left({\omega_{\varphi}}^{k+1} - {\omega}^{k+1}\right)\wedge
{\omega_{\varphi}}^{n-k-1}\wedge dt.
\]
Put $J_n =0 $ for convenience in notations.
\end{defi}

It is straightforward to verify that the definition of  $J_k(0\leq
k\leq n)$ is independent of path chosen.\\

Direct computations lead to
\begin{theo}\label{th:energydecay0} For any $k=0,1,\cdots, n$, we have
\begin{eqnarray}
{d E_k \over dt} & = & {{k+1}\over V} \displaystyle \int_M
\Delta_{\varphi}\left( {{\partial \varphi}\over {\partial t}}
\right )\; {{\rm Ric}(\omega_{\varphi})}^{k} \wedge
{\omega_{\varphi}}^{n-k}
\nonumber \\
& & \qquad -
 {{n-k}\over V}\displaystyle \int_M {{\partial \varphi}\over {\partial t}} \left({{\rm Ric}(\omega_{\varphi})}^{k+1}
 - {\omega_{\varphi}}^{k+1}\right) \wedge  {\omega_{\varphi}}^{n-k-1}.
\label{eq:decay functional0}
\end{eqnarray}
Here $\{\varphi(t)\}$ is any path in ${\cal P}(M,\omega)\;$ and
$V$ is the volume of K\"ahler metrics in this K\"ahler class
$[\omega]\;$
\end{theo}

\begin{prop} \label{th:energydecay}
Along the K\"ahler Ricci flow where $Ric > -1$ is preserved, we have
\begin{equation}
{d E_k \over d t} \leq - {{k+1}\over V} \displaystyle \int_M
(R(\omega_\varphi)-r) {\rm Ric}(\omega_{\varphi})^{k} \wedge
{\omega_{\varphi}}^{n-k}. \label{eq: Ekdecreases}
\end{equation}
When $k=0 , 1$, we have
\begin{eqnarray}
{{d E_0 }\over{d\,t}} & =  & -{{n\sqrt{-1}}\over V} \displaystyle
\int_M \partial {{\partial \varphi}\over {\partial t}} \wedge
\overline{\partial} {{\partial \varphi}\over{\partial
t}} \wedge  {\omega_{\varphi}}^{n-1}\leq 0,\label{eq:E0decrease}\\
{{d E_1 }\over{d t}} & \leq & - {{2}\over V} \displaystyle \int_M
(R(\omega_\varphi)-r)^2 {\omega_{\varphi}}^{n} \leq  0.
\label{eq:E1decrease}
\end{eqnarray}
In particular, both $E_0$ and $E_1$ are decreasing along the
K\"ahler Ricci flow.
\end{prop}

\section{Proof of Theorem 1}

\subsection{The average $L^2$ norm of the trace-less Ricci tensor}
\begin{lem}
Suppose that the curvature of $g_0$ satisfies the following
condition
\begin{equation}
\left\{\begin{array}{lcl} \mid R_{i\b j k \b l}(g_0) \mid & \leq & \Lambda > 0, \\
   R_{i \b j} (g_0) &\geq & - 1 + \delta > 0.\end{array}\right.
\end{equation}
Then, there exists a constant $T(\delta, \Lambda) > 0$, such that
the following bound for the evolved K\"ahler metric $g(t) (0\leq
t\leq 6 T)$
\begin{equation}
\left\{\begin{array}{lcl} \mid R_{i\b j k \b l}(g(t)) \mid & \leq & 2 \Lambda > 0, \\
   R_{i \b j} (g(t)) &\geq & - 1 + {\delta \over 2} > 0.\end{array}\right.
\end{equation}
In fact, we can choose $6T$ to be
\[
 \left( {1\over \Lambda} \right)^{2n} \delta \cdot \epsilon_0(n)^2.
\]
\end{lem}
\begin{proof} We use $BR$ to denote the bisectional
curvature, and ``*" to denote the multiplication and contraction
of two curvature tensors.  Then, the evolution equation
\ref{eq:evolutio of curvature1}  can be re-written schematically
\[
    {\p \over {\p t}} BR = \triangle BR + BR*BR
 \]
 If we let $u= \mid R_{i\b j k \b l}\mid^2$, then
 \[
  {\p \over {\p t}} u \leq \triangle u + c |Riem| \cdot  u = \triangle u + c \cdot u^{3\over 2} .
 \]
 where $c$ is some universal constants (depends on flow equation
 algebraically).
 This implies that
 \[
 \displaystyle \max_M \mid R_{i\b j k \b l}\mid \leq {1\over {\Lambda^{-1} - t}}.
 \]
 For $ 0\leq t \leq {1\over {2 \Lambda}}$, we have
 \[
    \displaystyle \max_M \mid R_{i\b j k \b l}\mid \leq 2 \Lambda.
 \]

Similarly, re-write the evolution equation \ref{eq:evolutio of
curvature2} for the Ricci curvature schematically:
\[{\p\over {\p t}} Ric = \triangle Ric + Ric* BR.\]
Applying the estimate of the bisectional curvature for $ 0\leq t
\leq 6 T,\;$ we have
 \[\left\{\begin{array}{lcl} \mid R_{i\b j k \b l}(g_0) \mid & \leq & 2 \Lambda > 0 \\
 R_{i \b j} (g_0) &\geq & -1 + {\delta\over 2} > 0.\end{array}\right.
  \]
 Here $c(n)$ is some universal constant and we can re-choose $6T$
 (may need to make it smaller if needed) as
  \begin{equation}
  6 T = {\delta\over {(\delta +\Lambda \cdot c(n))\cdot\Lambda}}.
  \label{eq:estimateontime}
 \end{equation}
\end{proof}

\begin{lem} \label{th:E1decay} If $$E_1(g_0) \leq  \displaystyle \inf_{Ric (\omega) \geq -1, \omega \in [\omega]} E_1
(g) + \epsilon(\delta, \Lambda), $$ then
\begin{equation}
{1\over {6 T}}  \displaystyle \int_0^{6T} \displaystyle \int_M\;
\mid Ric(\omega_\varphi) - \omega_\varphi\mid^2 d\,t \leq
{\epsilon_0(n)^2\over 2}.
\end{equation}
\end{lem}
\begin{proof}
For $0\leq t \leq 6T, $ we have $Ric_{\omega(g(t))} + \omega(g(t))
\geq {\delta\over 2} \geq 0.\;$  By  Theorem
\ref{th:energydecay0}, we have
\[
\begin{array}{lcl}
{d E_1 \over dt} & = & {2\over V} \displaystyle \int_M
\Delta_{\varphi}\left( {{\partial \varphi}\over {\partial t}}
\right )\; {{\rm Ric}(\omega_{\varphi})} \wedge
{\omega_{\varphi}}^{n-1}\\
& & \qquad -
 {{n-1}\over V}\displaystyle \int_M {{\partial \varphi}\over {\partial t}} \left({{\rm Ric}(\omega_{\varphi})}^{2}
 - {\omega_{\varphi}}^{2}\right) \wedge  {\omega_{\varphi}}^{n-2}\\
 & = & {2\over V} \displaystyle \int_M
(n - R(\omega(g(t)))) (R(\omega_{g(t)}) - n) \wedge
{\omega_{\varphi}}^{n-1}\\
& & \qquad -
 {{n-1}\over V}\displaystyle \int_M \sqrt{-1} \p {{\partial \varphi}\over {\partial t}} \wedge \b \p {{\partial \varphi}\over {\partial t}}
  \wedge \left({{\rm Ric}(\omega_{\varphi})}
 + {\omega_{\varphi}}\right) \wedge  {\omega_{\varphi}}^{n-2}
\\ & \leq &  - {2\over V} \displaystyle \int_M
(n - R(\omega_{g(t)}))^2  {\omega_{\varphi}}^{n} \leq 0.
\end{array}
\]

Consequently, we have
\[\begin{array}{lcl} 2 \displaystyle \int_0^{6T} \displaystyle \int_M\;
\mid Ric(\omega_\varphi) - \omega_\varphi\mid^2 \omega^n\; d\,t &
\leq &  - \displaystyle \int_0^{6T}\;  {{d\,E_1(g(t))}\over
{d\,t}} \;d\,t
\\
& \leq & E_1(g(0)) - E_1(g(6 T))  \\
&\leq & E_1(g(0)) - \displaystyle \inf_{\omega_g \in [g_0]}\; E(g) \\
& \leq & \epsilon(\Lambda, \delta).
\end{array}
\]
If
\[
 \epsilon(\Lambda, \delta) \leq  {1\over 2} \epsilon_0(n)^2 \cdot {\delta\over {(\delta +\Lambda  c(n))\Lambda}},
\]
and the definition of $6T$ (cf. equation \ref{eq:estimateontime}),
then
\[
{1\over {6 T}}  \displaystyle \int_0^{6T}\;d\,t  \displaystyle
\int_M\; \mid Ric(\omega_\varphi) - \omega_\varphi\mid^2
d\;vol_{g(t)} \leq {\epsilon_0(n)^2\over 2}.
\]
\end{proof}

\subsection{Estimates of Sobolev and Poincare constants}
The content of this subsection can be founded in
\cite{chentian002}. We re-produce it here for the convenience of
readers. In this section, we will prove that for any K\"ahler
metric in the canonical K\"ahler class, if the scalar curvature is
close enough to a constant in $L^2$ sense and if the Ricci
curvature is non-negative, then there exists a uniform upper bound
for both the Poincar\'e constant and the Sobolev constant. We
first follow an approach taken by  C. Sprouse \cite{spro001} to
obtain a uniform upper bound on the diameter.

 In \cite{CC96}, J. Cheeger
and T. Colding proved an interesting and useful inequality which
converts integral estimates along geodesic to integral estimates
on the whole manifold.  In this section, we assume $m= {\rm dim}
(M).\;$
\begin{lem}\label{th:lineintegral} \cite{CC96}
Let $A_1, \; A_2$ and $W$ be open subsets of $M$ such that $A_1,
A_2 \subset W,$ and all minimal geodesics $r_{x,y}$ from $x \in
A_1$ to $y \in A_2$ lie in $W.\;$ Let $f$ be any non-negative
function. Then
\[
\begin{array} {l} \displaystyle \int_{A_1 \times A_2} \displaystyle
\int_{r_{x,y}} \;  f(r(s)) \;d\,s\; d \; vol_{A_1 \times A_2}  \\
\qquad \leq C(m,k,\Re)({\rm diam}(A_2) vol(A_1) + {\rm diam}(A_1)
vol(A_2)) \displaystyle \int_W \; f\; d\;vol,
\end{array}
\]
where for $k \leq 0,$
\begin{equation}
C(m,k,\Re) = {{{\rm area} (\partial B_k (x,\Re))} \over {{\rm
area} (\partial B_k (x,{\Re \over 2}))}},
\end{equation}
\begin{equation}
\Re\geq \displaystyle \sup \{d(x,y) \mid (x,y) \in (A_1 \times
A_2)\},
\end{equation}
and $B_k(x,r)$ denotes the ball of radius $r$ in the simply
connected space of constant sectional curvature $k.\;$
\end{lem}
In this paper, we always assume ${\rm Ric} \geq 0$ on $M,$ and
thus $C(n,k,\Re) = C(n).\;$ \\
Using this theorem of Cheeger and Colding, C. Sprouse
\cite{spro001} proved an interesting lemma:

\begin{lem} \cite{spro001}
Let $(M,g)$ be a compact Riemannian manifold with ${\rm Ric} \geq
-1.\;$ Then for any $\delta > 0$ there exists $\epsilon =
\epsilon(n, \delta) $ such that if
\begin{equation}
{1 \over V} \displaystyle \int_M \left((m-1) - Ric_- \right)_+ <
\epsilon(m,\delta), \label{eq:dimaeter0}
\end{equation}
then the $diam(M) < \pi + \delta.$ Here ${\rm Ric}_-$ denotes the
lowest eigenvalue of the Ricci tensor; For any function $f$ on
$M,\; f_+(x) = \max\{f(x),0\}.$
\end{lem}

\begin{rem} Note that the right hand side of equation (\ref{eq:dimaeter0}) is not
scaling correct. A scaling correct version of this lemma should
be: For any positive integer $a> 0,$ if
\[
{1\over V}  \displaystyle \int_M \; |Ric- a| \; d\,vol <
\epsilon(m,\delta) \cdot a,
\]
then the diameter has a uniform upper bound.
\end{rem}

\begin{rem} It is interesting to see what the optimal constant
$\epsilon(m,\delta)$ is.  Following this idea, the best constant
should be
\[
\epsilon(m,\delta)  = \displaystyle \sup_{N > 2}\; {{N-2}\over {8
C(m) N^m}}.
\]
However, it will be  interesting to figure out the best constant
here. 
\end{rem}

Adopting his arguments, we will prove the similar lemma,
\begin{lem} \label{th:bestconstant2} Let $(M,\omega)$ be a polarized K\"ahler manifold and
$[\omega]$ is the canonical K\"ahler class. Then there exists a
positive constant $\epsilon_0(n)$ which only depends on the
dimension, such that if the Ricci curvature of $\omega$ is
non-negative and if
\[
 {1 \over V} \displaystyle \int_M \; (R -n)^2 \omega^n \leq
 \epsilon_0(n)^2,
\]
then there exists a uniform upper bound on diameter of the
K\"ahler metric $\omega.\;$ Here $r$ is the average of the scalar
curvature.
\end{lem}
\begin{proof}
We first  prove that the Ricci form is close to its K\"ahler form
in the $L^1$ sense (after proper rescaling).  Note that
\[ {\rm Ric}(\omega) - \omega = \sqrt{-1} \partial \bar \partial f
\]
for some real valued function $f.\;$ Thus
\[
\displaystyle \int_M\; \left( {\rm Ric}(\omega) - \omega \right)^2
\wedge \omega^{n-2} = \displaystyle \int_M\; \left( \sqrt{-1}
\partial \bar \partial f\right)^2 \wedge \omega^{n-2} =0.
\]

On the other hand,  we have
\[ \displaystyle \int_M\;
\left( {\rm Ric}(\omega) - \omega \right)^2 \wedge \omega^{n-2} =
{1\over {n(n-1)}}\; \displaystyle \int_M\; \left((R - n)^2 - |{\rm
Ric}(\omega) - \omega |^2 \right) \omega^n.
\]
Here we already use the identity $tr_{\omega}  \left( {\rm
Ric}(\omega) - \omega \right) = R - n.\;$ Thus
\[
\displaystyle \int_M\;|{\rm Ric}(\omega) - \omega |^2  \; \omega^n
=\displaystyle \int_M\;(R - n )^2  \; \omega^n.
\]
This implies that
\[
\begin{array}{lcl} \left(\displaystyle \int_M |{\rm Ric} - 1
|\;\omega^n\right)^2 & \leq & \displaystyle \int_M |{\rm
Ric}(\omega) -
\omega |^2 \;\omega^n \cdot \displaystyle \int_M\;\omega^n\\
& = & \displaystyle \int_M\; (R - n )^2 \;\omega^n \cdot V
\\ &\leq & \epsilon_0^2 \cdot V \cdot V = \epsilon_0^2 \cdot V^2,
\end{array}
\]
which gives
\begin{equation}
{1 \over V} \displaystyle \int_M |{\rm Ric} - 1 |\;\omega^n \leq
\epsilon_0. \label{eq:dimaeter2}
\end{equation}
The value of $\epsilon_0$ will be determined later.

Using this inequality (\ref{eq:dimaeter2}), we want to show that
the diameter must be bounded from above. Note that in our setting,
$m = {\rm dim}(M) = 2 n.\;$ Unlike in \cite{spro001}, we are not
interested in obtaining a sharp upper bound on the diameter. \\

Let $A_1$ and $A_2$ be two balla of small radius and $W = M. \;$
Let $f= |{\rm Ric} - 1| = \displaystyle \sum_{i=1}^m |\lambda_i -
1|,\;$ where $\lambda_i$ is the eigenvalue of the Ricci tensor. We
assume also that all geodesics are parameterized by arc length. By
possibly removing a set of measure $0$ in $A_1 \times   A_2,$
there is a unique minimal geodesic from $x$ to $y$ for all $(x,y)
\in A_1 \times A_2.$ Let $p,q$ be two points on $M$ such that
\[
   d(p,q) = {\rm diam}(M) = D.
\]

We also used $d\,vol$ to denote the volume element in the
Riemannian manifold $M$ and $V$ denote the total volume of $M.\;$
For $r> 0$, put $A_1 = B(p, r) $ and $A_2 = B(q,r).\;$ Then Lemma
\ref{th:lineintegral} implies that
\[
\begin{array} {l} \displaystyle \int_{A_1 \times A_2} \displaystyle
\int_{r_{x,y}} \; |{\rm Ric} - 1| \;d\,s\; d \; vol_{A_1 \times A_2} \nonumber \\
\qquad \qquad \leq C(n,k,R)(diam(A_2) vol(A_1) + diam(A_1)
vol(A_2)) \displaystyle \int_W \; |{\rm Ric} - 1| \;d\,vol.
\end{array}
\]
Taking infimum over both sides, we obtain
\begin{eqnarray} \displaystyle \inf_{(x,y) \in A_1 \times A_2 } \displaystyle
\int_{r_{x,y}} \; |{\rm Ric} - 1| \;d\,t\nonumber \\
\qquad \leq 2 \;r\; C(n)({1\over vol(A_1)} + {1\over  vol(A_2)})
\displaystyle \int_W \; |{\rm Ric} - 1|\;d\,vol
\nonumber \\
\leq 4 r C(n) {{D^n}\over r^n} {1\over V} \displaystyle \int_M \;
|{\rm Ric} - 1| \;d\,vol, \label{eq:CT2}
\end{eqnarray}
where the last inequality follows from the relative volume
comparison. We can then find a minimizing unit-speed geodesic
$\gamma$ from $x \in \overline{A_1}$ and $y \in \overline{A_2}$
which realizes the infimum, and will show that for $L =  d(x,y)$
much larger than $\pi,$ $\gamma$ can not be minimizing if the
right hand
side of (\ref{eq:CT2}) is small enough.\\

Let $E_1(t), E_2(t), \cdots E_m(t)$ be a parallel orthonormal
basis along the geodesic $\gamma$ such that $E_1(t) =
\gamma'(t).\;$ Set now $Y_i(t) = \sin \left({{\pi t}\over L}
\right) E_i(t), i =2,3,\cdots m.\;$ Denote by $L_i(s)$ the length
functional of a fixed endpoint variation of curves through
$\gamma$ with variational direction $Y_i,$ we have the 2nd
variation formula
\[
\begin{array}{lcl}
 & & \displaystyle
\sum_{i=2}^m {{d^2 L_i(s)}\over {d\,s^2}}\mid_{s=0} \\
& = & \displaystyle \sum_{i=2}^m \displaystyle \int_0^L \left(
g(\nabla_{\gamma'} Y_i, \nabla_{\gamma'}Y_i) -
R(\gamma',Y_i,\gamma',Y_i) \right)d\,t\\
& = & \displaystyle \int_0^L (m-1) \left({\pi^2 \over L^2} \cos^2
\left( {{\pi t}\over L} \right) \right)- \sin^2 \left( {{\pi
t}\over L}
\right) {\rm Ric}(\gamma',\gamma') \;d\,t \\
& = & \displaystyle \int_0^L \left((m-1) {\pi^2 \over L^2} \cos^2
\left( {{\pi t}\over L} \right) -  \sin^2 \left( {{\pi t}\over L}
\right) \right)\;d\,t \\
& & \qquad + \displaystyle \int_0^L\; \sin^2 \left( {{\pi t}\over
L} \right) \left( 1 - {\rm Ric}(\gamma',\gamma') \right)\;d\,t \\
& = & -{L\over 2} \left( 1- (m-1) {\pi^2 \over L^2} \right) \\
& & \qquad + \displaystyle \int_0^L\; \sin^2 \left( {{\pi t}\over
L} \right) \left( 1 - {\rm Ric}(\gamma',\gamma') \right)\;d\,t.
\end{array}
\]

Note that
\[
  1 -{\rm Ric}(\gamma',\gamma') \leq |{\rm Ric} - 1|.
\]

Combining the above calculation and the inequality (\ref{eq:CT2}),
we obtain

\begin{eqnarray} & & \displaystyle
\sum_{i=2}^n {{d^2 L_i(s)}\over {d\,s^2}}\mid_{s=0}  \nonumber \\
&\leq  & -{L\over 2} \left( 1- (m-1){\pi^2 \over L^2} \right)  +
\displaystyle \int_0^L\; \sin^2 \left( {{\pi t}\over L} \right)
|{\rm Ric} - 1|\;d\,t \nonumber\\
& \leq & -{L\over 2} \left(1- (m-1){\pi^2 \over L^2} \right)  + 4
r C(n) {{D^n}\over r^n} {1\over V} \displaystyle \int_M \; |{\rm
Ric} - 1| \;d\,vol. \label{eq:CT3}
\end{eqnarray}
Here in the last inequality, we have already used the fact that
$\gamma$ is a geodesic which realizes the infimum of the left side
of inequality (\ref{eq:CT2}).  For any fixed positive larger
number $N > 4$,  let $D= N\cdot r .\;$ Set $ c = {1 \over V}
\displaystyle \int_M \; |{\rm Ric} - 1| \;d\,vol.\;$  Note that
\textbf{}\[ L = d(x,y) \geq d(p,q) - 2 r = D (1- {2\over N})\geq
{D\over 2}.
\]
Then the above inequality (\ref{eq:CT3}) leads to
\[
\begin{array}{lcl}
  { 1 \over D}  \displaystyle
\sum_{i=2}^n {{d^2 L_i(s)}\over {d\,s^2}}\mid_{s=0}
 & \leq &  -{{1 - {2 \over N} }\over 2} \left(1- (m-1) {\pi^2 \over
L^2} \right) + 4 C(n) \; {N^{m-1} \over V} \cdot c  \cdot V\\
& = &  4 C(n) \; {N^{m-1}} \left( c - {{(N-2)}\over {2N}}{1 \over
{4 C(n) N^{m-1}}} \right)  + {{1 - {2 \over N} }\over
2}(m-1){\pi^2 \over L^2}.
\end{array}
\]
Note that the second term in the right hand side can be ignored if
$L \geq {D \over 2} $  is large enough.
Set
\[
\epsilon_0  = {{(N-2)}\over {2N}} \cdot {1 \over {4 C(n) N^{m-1}}}
= {{N-2 }\over {8 C(n) N^m  }}.
\]
Then if
\[
 {1 \over V} \displaystyle \int_M \; (R -n)^2 \omega^n \leq
 \epsilon_0^2,
\]
by the argument at the beginning of this proof, we have the
inequality (\ref{eq:dimaeter2}):
\[
 {1 \over V} \displaystyle \int_M \; (R -n)^2 \omega^n \leq
 \epsilon_0^2,
\] \[
{1\over V} \displaystyle \int_M \; |{\rm Ric} - 1| \;d\,vol <
\epsilon_0,
\]
which in turns imply
\[
 { 1 \over D}  \displaystyle
\sum_{i=2}^n {{d^2 L_i(s)}\over {d\,s^2}}\mid_{s=0} < 0,
 \]
for $D$ large enough. Thus, if the diameter is too
large,$\;\gamma$ cannot be a length minimizing geodesic. This
contradicts our earlier assumption that $\gamma$ is a
 minimizing geodesic. Therefore, the diameter must have a uniform upper bound.
\end{proof}

According to the work of.
 C. Croke \cite{croke80}, Li-Yau \cite{liyau80} and Li
\cite{pli80}), we state the following lemma on the upper bound of
the Sobolev constant and Poincare constant:
\begin{lem}\label{th:poincareconstant}
Let $(M,\omega)$ be any compact polarised K\"ahler manifold where
$[\omega]$ is the canonical class. If $Ric(\omega) \geq  - 1,\; V
= \displaystyle \int_M \omega^n  \geq \nu
> 0$ and the diameter has a uniform upper bound,
then there exists a constant $\sigma = \sigma(\epsilon_0, \nu)$
such that for all function $f\in C^{\infty}(M),\;$ we have
\[
\left( \displaystyle \int_M\;\mid f\mid^{{2n}\over {n-1}}\;
\omega^n \right)^{{n-1}\over n} \leq \sigma \left( \displaystyle
\int_M \mid \nabla f\mid^2 \;\omega^n + \displaystyle \int_M f^2
\;\omega^n \right).
\]
Furthermore, there exists a uniform Poincar\'e constant
$c(\epsilon_0) $ such that the Poincar\'e inequality holds
\[
\displaystyle \int_M \left( f - {1\over V} \displaystyle \int_M
\;f \;\omega^n \right)^2 \;\omega^n \leq c(\epsilon_0)\;
\displaystyle \int_M \; \mid\nabla f\mid^2 \; \omega^n.
\]
Here $\epsilon_0$ is the constant appeared in Lemma
\ref{th:bestconstant2}.
\end{lem}
\begin{proof} 
Note that $(M,\omega)$ has a uniform upper bound on the diameter.
Moreover,  it has a lower volume bound and it has non-negative
Ricci curvature. Following a proof in \cite{pli80} which is based
on a result of C. Croke \cite{croke80}, we obtain a uniform upper
bound on the Sobolev constant (independent of metric!).

Recall a theorem of Li-Yau \cite{liyau80} which gives a positive
lower bound of the first eigenvalue in terms of the diameter when
Ricci curvature is nonnegative: \[ \lambda_1(\omega) \geq
{\pi^2\over {4 D}}, 
\]
here $\lambda_1,\; D$ denote the first eigenvalue and the diameter
of the K\"ahler metric $\omega.\;$  In the same paper, Li-Yau also
have a lower bound control of the first eigenvalue in terms of the
Diameter when the Ricci curvature is bounded
from below.  \\

Now $D$ has a uniform upper bound according to Lemma
\ref{th:bestconstant2}. Thus the first eigenvalue of $\omega$ has
a uniform positive lower bound; which, in turn, implies that there
exists a uniform Poincar\'e constant.
\end{proof}

\begin{lem} Along the K\"ahler Ricci flow, the diameter of the evolving metric is uniformly bounded for $t\in [T,
4T].$
\end{lem}
\begin{proof} Lemma \ref{th:E1decay} implies implies that
\[
 {1\over {6T}} \displaystyle \int_0^{6T}\;d\,t \displaystyle \int_M\; \; (R-n) ^2\;
 \omega_{\varphi}^{n}\leq  {\epsilon_0(n)^2\over 2}.
\]
In particular, there exists at least one time $t_0 <  T$ such that
\[
 \displaystyle \int_M\; \; (R-n) ^2\;
 \omega_{\varphi(t_0)}^{n}\leq  {\epsilon_0(n)^2\over 2}.
\]

Now for this $t_0$,  applying lemma \ref{th:bestconstant2}, we
show there exists a uniform constant $D$ such that the diameters
of $\omega_{\varphi(t_0)}$ are uniformly bounded by ${D\over
2}.\;$ Recalled that $Ric(\omega) > -1$ for $t \in [0, 6T]\;$  so
that diameter of evolving metric increased at most exponentially
since
\[
  {\p\over {\p t}} g_{i\bar j} = g_{i \bar j} - R_{i \bar j} \leq
 2  g_{i\bar j}.
\]
Now $t - t_0 < 6 T,$ this implies that the Diameters of the
evolving metric along the entire flow is controlled by $e^{8T}
{D\over 2}.\;$
\end{proof}
\begin{rem} Notice that this diameter constant $D$ is not related
to $\delta, \Lambda.\;$ As $\epsilon(\delta, \Lambda) \rightarrow
0,$ we have $D\rightarrow \pi.\;$ This should be very useful
observation for future application.
\end{rem}
Combining this with Lemma \ref{th:poincareconstant}, we obtain

\begin{theo} Along the K\"ahler Ricci flow,  the evolving K\"ahler metric $\omega_{\varphi(t)} (t\in [T, 4T])$
has a uniform upper bound on the  Sobolev constant and Poincar\'e
constant.
\end{theo}
\subsection{The pointwise norm of traceless Ricci}
We first state a parabolic version of Moser iteration argument
\begin{lem} \label{lem:ricciiteration} If the Poincare constant
and the Sobolev constant of the evolving K\"ahler metrics $g(t)$
are both uniformly, and if a non-negative function $u$ satisfying
the following inequality
\[
{\p \over {\p t}} u \leq \triangle u + f(t,x) u, \qquad \forall a<
t < b,
\]
where $ f = \mid Riem(g(t))\mid$ and $|f|_{L^p(M,g(t))}$ is
uniformly bounded by some constant $c,$ then for any $\lambda \in
(0,1)$ fixed, we have
\[\displaystyle \max_{(1-\lambda) a + \lambda b \leq t \leq b} u \leq C(c, b-a, \lambda)\int_a^b\; \int_M u.
  \]
\end{lem}
\begin{rem} This is more or less standard and we refer readers to
literatures of evolution equations (c.f. Lemma 4.7
\cite{chentian002}).
\end{rem}
\begin{theo} \label{th:ricciiteration}  For any $\delta, \Lambda > 0, $ there exists a small positive constant $\epsilon (\delta, \Lambda) > 0$
 such that if the initial metric $g_0$ of the K\"ahler Ricci flow satisfies the following condition:
\begin{equation}
Ric(g_0) > -1 + \delta, |Riem(g_0)| < \Lambda, E_1(g_0) \leq \inf E_1 + \epsilon,
\label{eq:ricciIterationCond}
\end{equation}
 then there exists an interval length $T$ which depends on $\delta, \Lambda$ such that, after time $2T$ along the K\"ahler Ricci flow, we have
\[
\displaystyle\;\max_{t\in [2T, 6T]} \displaystyle\max_{x\in M} \;
\mid Ric_{g(t)} - \omega_{g(t)}\mid \leq \epsilon(n).
\]
\end{theo}
\begin{proof} Recalled that the evolution equation for Ricci
curvature is:
\[
 {{\p R_{i\b j}(g(t)) }\over {\p t}} = \triangle_{g(t)} R_{i\b
 j}(g(t)) + R_{p \b q}(g(t)) \cdot R_{q\b p i\b j}(g(t)) - R_{ i \b j}, \qquad
 \forall\; i, j =1,2,\cdots n.
\]
Set $Ric^0 = Ric -\omega. $ Then $u=\mid Ric^0\mid^2$ satisfies
this parabolic inequality
\[
{\p \over {\p t}} u \leq \triangle u + c_3 \mid Riem(g(t))\mid_g
\cdot u, \qquad \forall 0< t < 4T.
\]
where $c_3$ is some universal constant.  In this section, we use
$\mid Riem(g(t))\mid \leq 2 \Lambda$ for $0\leq t \leq 6T.\;$ Then
\[
\begin{array}{lcl} & & \displaystyle\;\max_{2T\leq t \leq 4 T} \displaystyle\;\max_{x\in M}\; \mid Ric^0(g(t))\mid \\
& = &
  \displaystyle \max_{2T\leq t \leq 4 T} \displaystyle \max_{x\in M}\;u(x,t) \\
  & \leq &  C(\Lambda,
4T,{1\over 4})\int_{T}^{4T}\;d\,t\; \displaystyle \int_M u
\;d\,vol_g \\ & = & C(\Lambda, 4T,{1\over 4})\int_{T}^{4T}\;d\,t\;
\displaystyle \int_M
|Ric(g(t))- g(t)|_g^2\;d\,vol_g \\
& \leq &  C(\Lambda, 4T,{1\over 4}) {\epsilon(n)_0^2\over 2} \cdot
4 T \leq \epsilon(n).
\end{array}
  \]
It is easy to see that we can extend this to proof the claim of
this lemma.
\end{proof}
In other words, the metric is nearly K\"ahler-Einstein at $t=
4T.\;$ If the norm of the bisectional curvature is still bounded
by $\Lambda, $ then we can continue to flow the metric from $4T$
to $ 8T$ to obtain a better estimate on ricci curvature, etc.
Unfortunately, we only have that the norm of the bisectional
curvature bounded from above by $2 \Lambda.\;$\\
\subsection{The proof of Theorem 1.}
We called a K\"ahler class satisfies condition (*) if it satisfies
the following condition:
\begin{equation}
[C_1(M)]^{[2]}\cdot [\omega]^{[n-2]} - {{2(n+1)}\over n} [ C_2(M)]
\cdot [\omega]^{[n-1]} = 0. \label{intro:topologycondition}
\end{equation}
 We actually only need the left hand side to be small enough for the present paper.  Note that the left side is always an integer and the only
manifold we know of which satisfies this condition is $CP^n.\;$
Consider its evolution equation:
\[
\begin{array}{lcl} {{\p R_{i\b j k\b l}}\over {\p t}} & =  & \triangle_{\varphi(t)} R_{i\b j k\b l} + + R_{i
\overline{j} p \overline{q}} R_{q \overline{p} k \overline{l}} -
R_{i \overline{p} k \overline{q}} R_{p \overline{j} q
\overline{l}} + R_{i \overline{l} p \overline{q}} R_{q
\overline{p} k \overline{j}} \\ & & \qquad \qquad -\left( R_{i\b
p} R_{p \overline{j} k \overline{l}} + R_{p\b j} R_{i \b p k\b l}
+ R_{k \b p} R_{i \b j p \b l} + R_{p \b l} R_{i \b j k\b p}
\right) - R_{i\b j k \b l}.
\end{array}
\]
Set
\[
 Q^0_{i\b j k \l} = R_{i \b j k \l} -{1\over {n+1}} (g_{i\b j} g_{k \b l} + g_{i \b l} g_{k \b
 j})
\]
Then, modulo $Ric^0 $ and $R-n$ (both can be controlled by $\mid
Q^0_{i \b j k\b l}\mid$), we have
\[
\begin{array}{lcl}
 & & \left( {\p \over {\p t}}- \triangle\right) Q_{i \bar j k \bar l}   \\
 & = &  Q_{i \overline{j} p \overline{q}} Q_{q \overline{p} k \bar
l} +Q_{i \overline{l} p \overline{q}} Q_{q \overline{p} k \bar j}
- Q_{i \overline{p} k \overline{q}} Q_{p \overline{j} q
\overline{l}}  \\ & & \qquad \qquad  - Q_{i \overline{j} k
\overline{l}}.
\end{array}
\]
Let $v = \mid Q^0_{i\b j k \l} \mid$. Then $v$ satisfies  this
parabolic inequality
\[
{\p \over {\p t}} v \leq \triangle v + \Lambda v, \qquad \forall
0< t < 4T.
\]
A calculation of E. Calabi \cite{calabi82} shows that in the
Canonical K\"ahler class which satisfies condition
\ref{intro:topologycondition}, we have
\[
\;\displaystyle \int_M \mid Q_{i\b j k\b l}\mid^2 =
\;\displaystyle \int_M \mid Ric^0\mid^2.
\]
Then
\[
\begin{array}{lcl} && \displaystyle\max_{3T\leq t \leq 4 T} \displaystyle\max_{x\in M}\; \mid Q_{i\b j k \b l}(g(t))\mid \\
& = &
  \displaystyle\max_{3T\leq t \leq 4 T} \displaystyle\max_{x\in M}\;v(x,t) \\
  & \leq &  C(\Lambda,
4T,2T)\int_{2T}^{4T}\; \displaystyle \int_M v(x,t) \leq C(\Lambda,
4T,2T)\int_{2T}^{4T}\; {1\over V}\;\displaystyle \int_M v^2\\ & =
&  (\Lambda, 4T,2T)\int_{2T}^{4T}\; {1\over V}\;\displaystyle
\int_M \mid Q_{i\b j k\b
l}\mid^2\\
& = &  (\Lambda, 4T,2T)\int_{2T}^{4T}\; {1\over V}\;\displaystyle
\int_M \mid Ric^0\mid^2
\\& \leq &  C(\Lambda, 4T,2T)\int_{2T}^{4T}\; {1\over V}\;
{\epsilon(n)_0^2\over 2} \cdot 4 T \leq \epsilon(n).
\end{array}
  \]

Consequently, we obtain
\begin{prop} The bisectional curvature has a pointwise pinching
estimate for $t \in [3T, 4T]$
\[
\displaystyle\max_{3T}^{4T}\; \displaystyle\;\max_{x\in M}\;\mid
R_{i \b j k \l} -{1\over {n+1}} (g_{i\b j} g_{k \b l} + g_{i \b l}
g_{k \b
 j})\mid < \epsilon(n)_0.
\]
\end{prop}

Combining Theorems 29 and Proposition 30, noticed that the
pinching condition for bisectional curvature, Ricci curvature as
well the $E_1$ norm all improved at time $t= 4T$ (comparing to
time $t=0$), repeatedly applying this procedure, we have
\begin{theo}  Along the K\"ahler Ricci flow, the bisectional
curvature becomes positive after finite time.  Moreover, there
exists a sequence of time $t_i \rightarrow \infty$ and $ t_{i+1} -
t_i \leq 1$ such that
\[ \begin{array}{lcl} 0 & = &
\displaystyle\lim_{i\rightarrow \infty}\;\displaystyle\;\max_{x\in
M}\;  \mid Ric_{g(t)}
-\omega_{g(t)}\mid_{t_i}  \\
& = & \displaystyle\lim_{i\rightarrow
\infty}\;\displaystyle\;\max_{x\in M}\;\mid R_{i \b j k \l}
-{1\over {n+1}} (g_{i\b j} g_{k \b l} + g_{i \b l} g_{k \b
 j})\mid_{t_i}. \end{array}
\]
\end{theo}

Appealing to a theorem in \cite{chentian001}, the flow converges
exponentially fast to a K\"ahler Einstein metric with constant
positive bisectional curvature.

\section{The proof of Theorem \ref{th:general} and \ref{th:2dim} }
\subsection{A compactness lemma}
In this section, we specialize in the complex surface (4 real
dimensional).  Consider the compactness of the following space
\[
\cA(c, C) = \{(M, g)\mid vol(M,g) = 1, (n-1) c g \leq Ric(g) \leq
(n-1) C g, \int_M \; |Riem(g)^2 \leq C \},
\]
here $c, C$ are two positive constants.  Then this is a compact
space in $W^{2,p}(M, g)$ space for any $p
> 1.\;$ In particular,

\begin{theo} \label{th:compactness} Suppose $X$ is a fixed holomorphic vector field in M.  Consider all the
K\"ahler metrics in $\cA(c,C)$ where $Im(X)$ is a Killing vector field.
For any $p > 1$, there exists a universal constant $\Lambda(n,
c,C, p)$ such that for all such K\"ahler metrics $g$, we have
\[
  \left( \int_M\; \mid Riem(g)\mid_g^p \;d\,vol_g \right)^{1\over
  p} \leq \Lambda(c, C, p).
\]
\end{theo}
\begin{proof}  For any metric $g\in \cA(c,C)$, we have $Ric(g) > (n-1) c.\;$ By Meyer's theorem, the diameter $Diam(g)$ of $g$ has a uniform upper bound
\[
   Diam(g) \leq {\pi\over \sqrt{c}}.
\]
Suppose the volume of Euclidean ball of radius $\rho$ is $\omega_n
\rho^n.\;$ By Bishop-Gromov relative volume comparison theorem,
for any $p\in M$, the volume ratio ${{vol(B_\rho(p))}\over {c_n
\rho^n}}$ is non-increasing function in $\rho>0.\;$ In particular,
this means the volume of geodesic ball of $g$ has a uniform upper
bound (in terms of Euclidean ball). On the other hand, we have
$vol(B_{Diam(g)})(p) = vol(M,g) = 1.\;$ In other words, we have
the volume ratio of geodesic ball has a uniform lower bound
\[
{c^{n\over 2} \over {\pi}^n} \leq {{vol(B_\rho(p))}\over {\omega_n
\rho^n}} \leq 1, \qquad \forall\; p\in M, \;\;{\rm and}\;\; 0\leq
\rho\leq Diam(g).
\]
The fact $vol(M,g) = 1$ also implies that $Diam(g)$ has a uniform
positive lower bound.
 These together with the uniform Ricci bound implies a
uniform positive lower bound of injectivity radius for metrics in
$\cA(c,C)\;$ which are invariant under action of the fixed Killing
vector field. Otherwise, there exists a sequence of
metrics $g_k \in \cA(c,C), \; k\in \cN$ such that the
injectivity radius of $g_k$ approaches to $0\;$ and $\cL_{Im(X)} g_k = 0.\;$
  We can re-scale
this sequence of metrics to be $\tilde{g_k}$ such that harmonic
radius is $1$ and the $L^\infty$ norm of Ricci curvature tends to
$0.\;$ Consequently, there is a subsequence of $\tilde{g_k}$ such
that it converges to some metric $g_\infty$ in a manifold
$M_\infty.\;$ The convergence in each harmonic ball of  radius $1$ is uniformly $C^{1,\alpha}.\;$
It is easy to see that $(M_\infty, g_\infty)$ is
complete K\"ahler manifold with flat Ricci curvature.  Moreover, there is a smooth holomorphic vector
field  in $M_\infty$ whose imaginary part acting isometrically on the metric $g_\infty.\;$ The  orbit of this isomeric action is either
all closed or all open.  In the case when the orbit is closed,
the length of each orbit is finite and the limit metric must admit a circle action. If the orbit of isometric action
is open, then each orbit must
have infinite length.  In the second case, since the $L^2$ norm of the Riemannian curvature of $g_\infty$
is uniformly bound,  the full Riemannian  curvature must vanish completely. Thus $g_\infty$ is
Euclidean and $M_\infty$  is either $R^4$ or $R^{3} \times
S^1.\;$ The cylinder case is ruled out because $g_\infty$
has maximum Euclidean volume growth. Consequently, the Harmonic radius of
the limit metric is infinity.
 Since the Harmonic radius is lower semi-continuous under $C^{1,\alpha}$ convergence,
the harmonic radius of $\tilde{g_k}$ for $k$ large enough must be
bigger  than $2.\;$ This is a contradiction since we assume initially that the
harmonic radius is $1$ for all $\tilde{g}_k, k \in \cN$! Consequently, the limit metric must admit an invariant circle action.  According to
Hitch antaze, in the global holomorphic coordinate of $M_\infty$, the limit metric $g_\infty$ can be expressed as
\[
   w^{-1}(d\,x^2 + d\, y^2) + w d\,z^2 + w \theta^2, \qquad w > 0.
\]
Here $\theta$ represents the connection 1-form of the circle action.  The fact of Ricci flat means
\[
 \triangle w = 0.
\]
Note that there is no positive harmonic function in $\cR^3,$ which follows that $w=constant > 0$ and the limit metric
is Euclidean again.  In particular, $M_\infty= \cC^2$ and the harmonic radius is again $\infty.\;$
As before, this contradicts with our earlier assumption that the harmonic radius of $\tilde{g}_k$ is only $1.\;$
Thus, there is a  uniform bound on the injectivity radius for this sequence of metrics $\{g_k, k \in \cN\}.\;$
Consequently,
 the
metric $g_k, k\in \cN$ is in $W^{2,p}(M, g)$ for any $p
> 1.\;$ In other words, we have
\[
 \left( \int_M\; \mid Riem(g)\mid_g^p \;d\,vol_g \right)^{1\over
  p} \leq \Lambda(c, C, p)
\]
for some uniform constant $\Lambda(c,C, p).$
 \end{proof}

\subsection{Moser iteration in the ``move"}
In this subsection, we want to use Theorem \ref{th:ricciiteration}
repeatedly to obtain global pinching of Ricci tensor over time
$t\in [T,\infty).\;$ By Theorem \ref{th:ricciiteration}, we have
\[
   1-\epsilon_n < Ric(g(t)) \leq 1+\epsilon_n, \qquad \forall t\in
   [2T,6T].
\]
\begin{lem} \label{le:mu-iteration} Fix a constant $p > {{2n+2}\over 2} = n+1.\; $ Condition as in Thereom \ref{th:ricciiteration}. Along the K\"ahler Ricci flow, if $A> 6T$ is any
positive number such that
\begin{equation}
   1-\epsilon_n < Ric(g(t)) \leq 1+\epsilon_n, \qquad \forall t\in
   [2T,A]
\label{eq:mu-iteration1}
\end{equation}
and
\begin{equation}
 \int_M\; |Riem(g(t))|^p \leq \Lambda({1\over 2}, 2, p),\qquad \forall t\in
   [2T,A],
\label{eq:mu-iteration2}
\end{equation}
then there exists a $\epsilon_1 > 0, $ such that both
inequalities hold for $t\in [2T, A+\epsilon_1].\;$
\end{lem}
\begin{proof}

Because that the K\"ahler Ricci flow exists globally and the metrics evolved smoothly , there is a
small constant $\epsilon_1$ such that
\[
 0 < (n-1)c < Ric(g(t)) \leq (n-1) C, \qquad \forall \;t\in [2T,
  A+\epsilon_1].
\]

Note that $Ric(g(t)) > -1$ for any $t \leq A+ \epsilon_1.\;$ Thus,
$E_1$ is decreasing and
\[
\begin{array}{lcl}
\displaystyle \int_0^{A+\epsilon_1} \;d\,t\displaystyle \int_M\;
\mid Ric(\omega_{g(t)}) -\omega_{g(t)}\mid^2_{g(t)}
\;d\,vol_{g(t)}&  \leq & E_1(0) - E_1(A+\epsilon_1)\\
&  \leq & E_1(0) - \displaystyle \inf_{[\omega]} \; E_1(g) <
\epsilon(\delta, \Lambda).
\end{array}
\]
On the other hand,  applying the compactness Theorem
\ref{th:compactness} to metric $g(t) (2T\leq t \leq
A+\epsilon_1)$, we have
\[
  \mid f(x,t)\mid_{L^{p}(M,g(t))} \leq \Lambda({1\over 2},2, p), \qquad {\rm here}\; f(x,t) = \mid Riem(g(t))\mid.
\]
 For $t\in [A-4T
+\epsilon_1, A+\epsilon_1],$ applying Lemma \ref{lem:ricciiteration}, we obtain
\[
\displaystyle\;\max_{t\in [A- 4T+\epsilon_1, A+\epsilon_1]}
\displaystyle\max_{x\in M} \; \mid Ric_{g(t)} -
\omega_{g(t)}\mid_g \leq \epsilon(n).
\]
Here we already use the fact that
\[
\begin{array}{lcl}
  \displaystyle \int_{A- 4T+\epsilon_1}^{A+\epsilon_1}
  \;d\,t\;\displaystyle\int_M\;\mid Ric_{g(t)} -
  \omega_{g(t)}\mid_{g(t)}^2 \;d\,vol_g & \leq & E_1(A- 4T+\epsilon_1) -
  E_1(A +\epsilon_1) \\
  & < & E_1(0) - \displaystyle \inf_{[\omega]} E_1(g) \\
  & \leq & \epsilon(\delta, \Lambda) \leq {\epsilon_0(n)^2\over 2}\cdot 5 T.
\end{array}\]

\end{proof}
Hence, we can repeatedly applied this iteration Lemma to itself.  This Lemma essentially says that the time $t$ where both
 conditions \ref{eq:mu-iteration1}
and \ref{eq:mu-iteration2} to hold is relatively open in $[2T, \infty).\;$  It is easy to show that it is relatively
closed as well.   Therefore,  we have the following direct corollary.
\begin{cor} Along the K\"ahler Ricci flow, we have
\[
\displaystyle\;\max_{t\in [2T,\infty)} \displaystyle\max_{x\in M}
\; \mid Ric_{g(t)} - \omega_{g(t)}\mid \leq \epsilon(t) \rightarrow 0.
\]
Moreover, we have a uniform bound of $Riem(g)$ for
any $t\in [0,\infty).\;$
\end{cor}
\begin{proof}  It is easy to see that the pinching condition for Ricci curvature holds. Using the compactness theorem \ref{th:compactness}, we show that $L^p$ norm of the Riemannian
curvature is uniformly bounded over the entire flow.  Using the iteration scheme as in section 3.4, but working on a longer interval, we can
show that the full Riemannian curvature is uniformly bounded.
\end{proof}
Now we are ready to wrap up the proof to Theorem  \ref{th:general} and \ref{th:2dim} .
\begin{proof}
According to the preceding corollary,  the Riemannian curvature is uniformly bounded
and the Ricci curvature uniformly pinched towards a  K\"ahler Einstein metric. Therefore,
for any sequence $t_i\rightarrow \infty,$ there exists a subsequence such that $g(t_i)$
converges to a K\"ahler Einstein metric, perhaps in a different
complex manifold. However, the stability assumption in general dimension will enable us
to derive the exponential convergence of the flow, similar to what is done in
\cite{chentian001}. This completes the proof of Theorem \ref{th:general}. In complex surface, we can bypass this stability conditions because
the classification theorem: Only K\"ahler surface with positive first Chern class are $CP^2, CP^1\times CP^1$
and $CP^2$ blowup $k (3 \leq k \leq 8)$ points in generic position (means no three points in one line
and no six points in a quadratic).  To have a non-trivial reductive automorphism group, we need further restricted
to $k = 3,4 $ or $CP^2, CP^1 \times CP^1.\;$ To have a Killing vector field, then the case $k=4$ is ruled out. In either of the remaining cases,
the orbit of complex structure under the
diffeomorphism group is well known.  We can show that the complex structure is stable
in the sense described in Theorem \ref{th:general}. This complete the proof for Theorem \ref{th:2dim}.
\end{proof}

\section{Future problems}
The following problems are interesting and will be investigated in the sequel of this paper.
\begin{enumerate}
\item Can we prove a general convergence theorem without assuming small energy on $E_1$? This is still interesting even if we assume
the condition \ref{intro:topologycondition}.
\item In K\"ahler Einstein manifold, is $E_1$ necessary has a lower bound or proper? Is the critical point of $E_1$ unique
(up to holomorphic transfermation)?  In other words, is the critical points necessary K\"ahler Einstein?
\item Under what geometric condition, can we derive the existence of lower bound of $E_1$ in absence of K\"ahler Einstein metric?
In particular, what about the special case where the bisectional curvature is positive?
\end{enumerate}

Xiuxiong Chen, Department of Mathematics/ University of Wisconsin/
Madison, WI 53706/ USA/ xiu@math.wisc.edu
\end{document}